\documentclass[10pt,intlimits]{amsart}

\usepackage[hypertex]{hyperref}
\usepackage{amsfonts,amssymb}
\usepackage{amsthm}
\usepackage{bbm} 
\usepackage{txfonts}

\usepackage{mathrsfs}

\usepackage{a4wide}

\newtheorem{neu}{}[section]

\newtheorem*{Cor*}{Corollary}

\newtheorem*{Thm*}{Theorem}
\newtheorem{Prop}[neu]{Proposition}
\newtheorem*{Prop*}{Proposition}
\theoremstyle{definition}

\newtheorem*{Rmk*}{Remark}
\newtheorem{Rmk}[neu]{Remark}

\newtheorem*{Ex*}{Example}

\newtheorem{Def}[neu]{Definition}

\newcommand{\Z}{\mathbb{Z}}
\newcommand{\R}{\mathbb{R}}
\newcommand{\C}{\mathbb{C}}

\newcommand{\pf}{\longrightarrow}

\newcommand{\wrt}{with respect to }

\newcommand{\Mas}{\mu_{\mathrm{Maslov}}}

\newcommand{\im}{\mathrm{image\,}}

\newcommand{\ev}{\mathrm{ev}}

\newcommand{\D}{\mathbb{D}}

\newcommand{\M}{\mathcal{M}}
\newcommand{\Mh}{\widehat{\mathcal{M}}}

\renewcommand{\H}{\mathrm{H}}

\newcommand{\PS}{\mathcal{PS}(S)}
\newcommand{\PSh}{\widehat{\mathcal{PS}}(S)}

\newcommand{\beq}{\begin{equation}}
\newcommand{\beqn}{\begin{equation}\nonumber}
\newcommand{\eeq}{\end{equation}}

\newcommand{\bea}{\begin{equation}\begin{aligned}}
\newcommand{\bean}{\begin{equation}\begin{aligned}\nonumber}
\newcommand{\eea}{\end{aligned}\end{equation}}

\numberwithin{equation}{section}

\begin{document}
\title{On the Weinstein conjecture in higher dimensions}
\author{Peter Albers}
\author{Helmut Hofer}
\thanks{Research partially supported by NSF grant DMS-0603957 and DFG grant AL 904/1-1}
\address{
    Peter Albers\\
    Courant Institute\\
    New York University}
\email{albers@cims.nyu.edu}
\address{
    Helmut Hofer\\
    Courant Institute\\
    New York University}
\email{hofer@cims.nyu.edu}
\begin{abstract}
The existence of a ``Plastikstufe'' for a contact structure implies the Weinstein conjecture for all supporting contact forms.
\end{abstract}
\maketitle

\section{Introduction and Main Result}
A one-form $\lambda$ on an odd-dimensional manifold $M^{2n-1}$ is
called a contact form, provided $\lambda\wedge d\lambda^{n-1}$ is a
volume-form. Associated to a contact form $\lambda$ we have the Reeb
vector field $X$ defined by
$$
i_X\lambda=1\ \ \hbox{and}\ \  i_Xd\lambda=0
$$
and the contact  structure $\xi=\hbox{ker}(\lambda)$. In 1978, A.
Weinstein, \cite{Weinstein_The_conjecture}, motivated by a result of P. Rabinowitz, \cite{Rabinowitz},
and one of his own results, \cite{Weinstein_Annals}, made the following
conjecture:\\

\begin{center} {\bf
A Reeb vector field on a closed manifold $\boldsymbol{M^{2n-1}}$
admits a periodic orbit.}
\end{center}
\vspace{0.3cm} The first break-through on this conjecture was
obtained by C.~Viterbo, \cite{Viterbo_Proof_of_Weinstein}, showing that compact energy
surfaces in ${\mathbb R}^{2n}$ of contact-type have periodic orbits.
Extending Gromov's theory of pseudoholomorphic curves, \cite{Gromov_Pseudoholomorphic_curves_in_symplectic_manifolds}, to
symplectized contact manifolds, H.~Hofer, \cite{Hofer_Weinstein_conjecture}, related the
Weinstein conjecture to the existence of certain pseudoholomorphic
curves. He showed that in dimension three the Weinstein conjectures
holds in many cases. In particular, he showed
that Reeb vector fields associated to over-twisted contact
structures admit periodic orbits. Recently the Weinstein conjecture
in dimension three was completely settled by C.~Taubes, \cite{Taubes_Weinstein_conjecture_I,Taubes_Weinstein_conjecture_II}, who
exploited relationships between Seiberg-Witten-Floer homology,
\cite{Kronheimer_Mrowka_book}, and embedded contact homology, \cite{Hutchings_index_inequality}, in order to
construct holomorphic curves in the symplectized contact manifold
out of nontrivial Seiberg-Witten-Floer homology classes. For more references on the Weinstein conjecture see \cite{HWZ_Dynamics_in_dim_three}.

In this note we show that many Reeb vector fields on higher
dimensional closed manifolds have periodic orbits generalizing  the
main result from \cite{Hofer_Weinstein_conjecture}. Our existence result
is closely connected to the interesting attempt by K.~Niederkr\"uger \cite{Niederkruger_Plastikstufe} to
generalize the three-dimensional notion of an overtwisted contact
structure. He introduced the concept of a
\textit{Plastikstufe} which currently seems to be the most
compelling generalisation given recent further developments by F.~Presas,
\cite{Presas_A_class_of_non_fillable_contact_structures} and K.~Niederkr\"uger / O.~van Koert, \cite{Niederkruger_van_Koert_Every_contact_manifold}.

Let us denote by $(M,\xi)$ a pair consisting of a closed manifold
$M$ of dimension $2n-1$ and a co-oriented  contact structure $\xi$.
We denote by $\D^2$ the closed unit disk in ${\mathbb C}$ with
coordinates $x+iy$.
\begin{Def}
We say that $(M,\xi)$ contains a \textit{Plastikstufe} with singular
set $S$ provided $M$ admits a closed submanifold $S$ of dimension
$n-2$ and an embedding $\iota:\D^2\times S\rightarrow M$ with
$\iota(\{0\}\times S)=S$ having the following properties:
\begin{enumerate}
\item There exists a contact form $\lambda_{PS}$ inducing $\xi$ so that the one-form
$\beta:=\iota^\ast\lambda_{PS}$ satisfies $\beta\wedge d\beta=0$ and
moreover $\beta\neq 0$ on $(\D^2\setminus\{0\})\times S$. Near
$\{0\}\times S$ the form $\beta$ is given by $\beta=xdy-ydx$ and the
pull-back of $\beta$ to $\partial\D^2\times S$ vanishes.
\item The complement of $\{0\}\times S$ in $(\D^2\setminus\partial\D^2)\times S$
 is smoothly foliated by $\beta$ via an $S^1$-family of leaves diffeomorphic
 to $(0,1)\times S$, where one of the ends converges to the singular set
 $\{0\}\times S$ and the other is asymptotic to the leave $\partial\D^2\times
 S$.
\end{enumerate}
The set $\PS=\iota(\D^2\times S)$ is called the
\textit{Plastistufe}.
\end{Def}
Let us observe that the existence of a \textit{Plastikstufe} for a
given contact structure involves the existence of a certain inducing
contact form. This is different from the three-dimensional case
where an over-twisted disk is defined only in terms of the contact
structure and does not require the existence of a particular contact
form. In the following we shall call a closed co-oriented contact
manifold $(M,\xi)$ PS-overtwisted provided there exists a contact
form $\lambda_{PS}$ inducing $\xi$ containing a \textit{Plastistufe}.
Recently Niederkr\"uger and van Koert showed that every
odd-dimensional sphere $S^{2n-1}$ with $n\geq 3$ has a contact
structure admitting a \textit{Plastikstufe}. If now $(M^{2n-1},\xi)$
is a co-oriented contact manifold then a connected sum with an
PS-overtwisted sphere admits by standard arguments a contact
structure which is PS-overtwisted. In particular, any closed manifold
of dimension $2n-1$ admitting a co-oriented contact structure also
admits a PS-overtwisted contact structure. Our main result is the
following theorem.
\begin{Thm*}
Let $(M,\xi)$ be a closed PS-overtwisted contact manifold. Then
every Reeb vector field associated to a  contact form $\lambda$
inducing $\xi$ has a contractible periodic orbit.
\end{Thm*}

\begin{Rmk}
In \cite{Niederkruger_Plastikstufe} Niederkr\"uger shows that a PS-overtwisted contact structure does not have a semi-positive symplectic
filling. We noticed that some of his idea combined with ideas from \cite{Hofer_Weinstein_conjecture} lead to the above theorem.
We also observed that the limitation to semi-positive fillings is not necessary and can be removed using polyfolds \cite{Hofer_CMD_2004_Harvard}.
This will be discussed in a forthcoming paper.
\end{Rmk}

\section{Background}

\noindent All material in this section is taken from \cite{Niederkruger_Plastikstufe}.

\subsection{Local normal form}

Let $(M,\lambda)$ contain a Plastikstufe $\PS$. In \cite[section 3.1]{Niederkruger_Plastikstufe} it is proved that there exist
constants $\varepsilon,C>0$ and an open set $V$ in the symplectic manifold $\big((-\varepsilon,0]\times M,d(e^s\lambda)\big)$ such
that $\{0\}\times S\subset V$ and $V$ is symplectomorphic to the set
\beq
U:=\left\{\;\big(\:\!(z_1,z_2),(q,p)\:\!\big)\in\C^2\times T^*S\left|\begin{aligned}\;& -C<\mathrm{Re}(z_1)\leq0,\;-C<\mathrm{Im}(z_1)<C\\
&\;\mathrm{Re}(z_1)+\tfrac14|z_2|^2+\tfrac12||p||^2\leq0\end{aligned}\right.\;\right\}\,,
\eeq
in $\C^2\times T^*S$ which carries its natural symplectic structure. Moreover, $M\cap V$ corresponds to equality in the last equation and $\PS\cap V$
to equality and $\mathrm{Im}(z_1)=0$, $p=0$.

\subsection{Bishop family}

The local model $U$ contains a natural $(n-1)$-dimensional \textit{Bishop family} given by
\bea
u_{t_0,q_0}:\D^2&\pf \C^2\times T^*S\\
z&\mapsto\big((-t_0,2\sqrt{t_0}z),(q_0,0)\big)
\eea
where $0\leq t_0<C$ is a real parameter and $q_0\in S$. The maps $u_{(t_0,q_0)}$ are $(i\times j)$-holomorphic, where $j$ denotes the natural almost
complex structure on $T^*S$ induced by the Levi-Civita connection of a Riemannian metric on $S$. Moreover, they have boundary on the set corresponding
to $\PS$.

We denote by $J$ the almost complex structure on $V$ obtained by pulling back the almost complex structure $i\times j$ from $\C^2\times T^*S$. Then
we can pull back the Bishop family to holomorphic maps (denoted by the same symbols)
\beq
u_{t_0,q_0}:\D^2\pf V\subset(-\varepsilon,0]\times M\,.
\eeq

\subsection{Uniqueness results for holomorphic disks}

We extend the almost complex structure $J$ from the set $V$ to a compatible almost complex structure on $\big(W:=(-\infty,0]\times M,d(e^s\lambda)\big)$.
We introduce the following notation
\beq
\PSh=\PS\setminus(\partial\PS\cup S)
\eeq
and remark that $\PSh$ is totally real \wrt $J$. The following proposition is taken from \cite[Proposition 7]{Niederkruger_Plastikstufe}.

\begin{Prop}\label{prop:local_uniqueness}
Let $u:(\D^2,\partial\D^2)\pf(W,\PSh)$ be a $J$-holomorphic disk which is simple. Moreover, we assume that $u(S^1)\subset\PS$ bounds a disk in $\PS$ and
\beq
\im(u)\cap V\neq\emptyset\,.
\eeq
Then, up to an element in $\mathrm{Aut}(\D^2)$, we have
\beq
u=u_{t_0,q_0}\,,
\eeq
that is, after reparametrization, the holomorphic disk $u$ is a member of the Bishop family.
\end{Prop}

\section{Proof of the theorem}

\noindent By assumption there exists a contact form $\lambda_{PS}$ on $M$ containing a Plastikstufe. Let $\lambda$ be another contact form
inducing the same contact structure.

\begin{center}\it We assume by contradiction that there exists no contractible closed Reeb orbit for $\lambda$.\end{center}

\subsection{The set-up}

We choose a function $f:M\pf\R$ such that $\lambda=f\lambda_{PS}$.
Since multiplying $\lambda$ with a non-zero constant doesn't change its Reeb orbits (up to reparametrization) we may assume without
loss of generality that the function $f$ takes only values in $(0,1)$. Then we can choose a smooth family of functions $f_s:M\pf\R$
for $s\in[-1,-\varepsilon]$ satisfying
\beq
f_s=\begin{cases}1&\text{near }s=-\varepsilon\\f&\text{near }s=-1\end{cases}\qquad\text{and moreover}\qquad\frac{\partial f_s}{\partial s}\geq0\;.
\eeq
This gives rise to a smooth family $\lambda_s=f_s\lambda_{PS}$ of contact forms which we extend by $\lambda_{PS}$ for $s\geq-\varepsilon$ and by $\lambda$
for $s\leq-1$. On $W=(-\infty,0]\times M$ we choose an exact symplectic form $\Omega$ on $W$ which satisfies
\beq
\Omega=\begin{cases}d(e^s\lambda_{PS})&\text{on }[-\varepsilon,0]\times M\\d(e^s\lambda)&\text{on }(-\infty,-1]\times M\end{cases}
\eeq
This is possible due to the choice of the family $f_s$. This has been used in the literature many times, see for instance
\cite{HWZ_Finite_energy_foliations_of_tight_three_spheres}.
We modify the almost complex structure $J$ from above to a compatible almost complex
structure $J$ on $(W,\Omega)$ by requiring that on $(-\infty,-2]\times M$ the almost complex structure is adapted to the negative part of the
symplectization of $\lambda$, in the sense of \cite{BEHWZ_SFT_compactness}. On $V$ it remains as defined in the
previous section. In particular, $(W,\Omega)$ still contains the Bishop
family $u_{t_0,q_0}$. We denote the relative homotopy class given by the Bishop disks by $a\in\pi_2(W,\PSh)$ and set
\begin{align}
\M(J)&:=\{u:(\D^2,\partial\D^2)\pf (W,\PSh)\mid\bar{\partial}_Ju=0,\,[u]=a,\,\mathrm{lk}(u,S)=1\}\,,\\
\Mh(J)&:=\M(J)\big/\mathrm{Aut}(\D^2)
\end{align}
where $\mathrm{lk}(u,S)$ is the linking number of $u(S^1)$ in $\PS$ with the set $S$. This is defined as follows.
By definition $\PSh$ is foliated by an $S^1$-family of Legendrian submanifolds,
thus there exists a natural map $\theta:\PSh\pf S^1$. We set $\mathrm{lk}(u,S):=\deg (\theta\circ u|_{S^1})$.

\subsection{The proof}

We need the following three facts established in \cite[Propositions 8 -- 10]{Niederkruger_Plastikstufe}.
\begin{enumerate}
\item The Maslov index of $a$ equals $\Mas(a)=2$,
\item the almost complex structure $J$ is regular at members of the Bishop family,
\item the energy of all elements in $\M(J)$ is uniformly bounded.
\end{enumerate}
The totally real submanifold $\PSh$ is non-compact. Since $\partial\PS$ is a closed leaf of the characteristic foliation the maximum principle
implies that no holomorphic maps intersect $\partial\PS$ at an interior point. According to Proposition \ref{prop:local_uniqueness}
near $S$ the only holomorphic disks are members of the Bishop family. Therefore, the non-compactness of $\PSh$ poses no problem. Moreover, due to
the energy bounds and the specific structure of the almost complex structure $J$ on the end of $W$ we can apply the ideas of the SFT-compactness theorem
\cite{BEHWZ_SFT_compactness}. Since we assumed that there exists no contractible closed Reeb orbits bubbling-off cannot occur in the interior.
Therefore, the only non-compactness of the moduli space $\Mh(J)$ comes from bubbling-off of holomorphic disks having boundary on $\PSh$.
The next proposition is taken from
\cite[Proposition 11]{Niederkruger_Plastikstufe} and shows that there exists no bubbling-off of holomorphic disks.
\begin{Prop}
Given a sequence $(u_n)\subset\Mh(J)$ there exists a subsequence either converging to an element in $\Mh(J)$ or to a point in $S$.
\end{Prop}
The latter case occurs if a family of Bishop disks shrinks to a point in $S$. We remark that in the former case the limit is simple.
\begin{Prop}
For a compatible almost complex structure $J$, which is generic on the subset $\big((-2,0]\times M\big)\setminus V$ of $(W,\Omega)$, the moduli
space $\M(J)$ is a smooth, compact manifolds of dimension
\beq
\dim\M(J)=n+2\,.
\eeq
\end{Prop}

\begin{proof}
We pick $u\in\M(J)$. In case that $\im(u)\cap V\neq\emptyset$ we conclude from Proposition \ref{prop:local_uniqueness} that
$u$ is a member of the Bishop family. In particular, $\im(u)\subset V$. Moreover, $J$ is already regular for members in the Bishop family.

If $\im(u)\cap V=\emptyset$ then it has to pass through the region $\big((-2,0]\times M\big)\setminus V$. Since all the disks are simple
a generic $J$ will be regular, see for example \cite{Dragnev}. The dimension formula follows from the fact that $\Mas(a)=2$ and $\dim\PSh=n$.
\end{proof}
We consider the evaluation map
\bea
\ev:\Mh(J)_{S^1}:=\M(J)\times_{\mathrm{Aut}(\D^2)}S^1&\pf\PSh\subset M\\
[u,t]&\mapsto u(e^{2\pi it})
\eea
defined on the smooth manifolds $\Mh(J)_{S^1}$ of dimension $\dim\Mh(J)_{S^1}=n$.
\begin{Prop}
For a generic $J$ as in the previous proposition the evaluation map is smooth.
\end{Prop}
To derive the contradiction to the assumption that $\lambda$ has no closed Reeb orbits we make the following
\begin{Def}
For a point $p=\iota(z,s)\in\iota(\D^2\times S)=\PS$ we define the distance of $p$ to $S$ by $d(p,S)=|z|$ and set for $0<\delta<\varepsilon$
\beq
\Mh(J)^{\delta}_{S^1}:=\Big\{[u,t]\in\Mh(J)_{S^1}\mid d(\ev([u,t]),S)\geq\delta\Big\}\;.
\eeq
\end{Def}
Then we have
\beq
\ev\big(\partial\Mh(J)^{\delta}_{S^1}\big)=\iota(S^1_\delta\times S)\,,
\eeq
where $S^1_\delta=\{z\in\D^2\mid|z|=\delta\}$. We conclude that $[\ev\big(\partial\Mh(J)^{\delta}_{S^1}\big)]\in\H_{n-1}(\PSh,\Z/2)$ is the
generator. On the other hand the set $\ev\big(\partial\Mh(J)^{\delta}_{S^1}\big)$ is clearly the boundary of the compact manifold
$\ev\big(\Mh(J)^\delta_{S^1}\big)$. This implies, that $[\ev(\partial\Mh(J)^{\delta}_{S^1})]=0\in\H_{n-1}(\PSh,\Z/2)$.

This contradictions concludes the proof of the theorem.

\noindent\hrulefill

%
%
\bibliographystyle{amsplain}

\bibliography{../../../Bibtex/bibtex_paper_list}
\end{document}